\documentclass[twoside,accepted]{article}
\usepackage{aistats2012}
\usepackage{natbib}
\usepackage{amsmath}
\usepackage{amssymb}
\usepackage{amsthm}
\usepackage{amsfonts}
\usepackage{latexsym}
\usepackage{epsfig}
\usepackage{times}
\usepackage{subfigure}

\newcommand{\rmPi}{\mathrm{\Pi}}
\newcommand{\iid}{\stackrel{iid}{\sim}}
\newcommand{\m}{\boldsymbol{\mathrm{m}}}
\newcommand{\Y}{\boldsymbol{\mathrm{Y}}}

\begin{document}

\twocolumn[

\aistatstitle{Stick-Breaking Beta Processes and the Poisson Process}

\aistatsauthor{John Paisley$^1$\quad David M. Blei$^3$\quad Michael I. Jordan$^{1,2}$}
\runningauthor{John Paisley\quad David M. Blei\quad Michael I. Jordan}

\aistatsaddress{$^1$Department of EECS, $^2$Department of Statistics, UC Berkeley\\$^3$Computer Science Department, Princeton University} ]

\begin{abstract}
We show that the stick-breaking construction of the beta process due to
\cite{Paisley:2010} can be obtained from the characterization
of the beta process as a Poisson process.  Specifically, we show that 
the mean measure of the underlying Poisson process is equal to that of the beta process. 
We use this underlying representation to derive error bounds on truncated beta processes
that are tighter than those in the literature. We also develop a new MCMC inference algorithm for beta processes, based in part on our new Poisson process construction.
\end{abstract}

\section{Introduction}
The beta process is a Bayesian nonparametric prior for sparse collections
of binary features~\citep{Thibaux:2007}.  When the beta process is marginalized 
out, one obtains the Indian buffet process (IBP)~\citep{Griffiths:2006}. Many applications of this circle of ideas---including focused 
topic distributions~\citep{Williamson:2010}, featural representations of 
multiple time series~\citep{Fox:2010} and dictionary learning for image 
processing~\citep{Zhou:2011}---are motivated from the IBP representation. However,  
as in the case of the Dirichlet process, where the Chinese restaurant 
process provides the marginalized representation, it can be useful to develop
inference methods that use the underlying beta process.  A step in this
direction was provided by \cite{Teh:2007a}, who derived a stick-breaking 
construction for the special case of the beta process that marginalizes
to the one-parameter IBP. 

Recently, a stick-breaking construction of the full beta process was 
derived by \cite{Paisley:2010}. The derivation relied on a limiting
process involving finite matrices, similar to the limiting process used to derive the IBP.
However, the beta process also has an underlying
Poisson process~\citep{Jordan:2010,Thibaux:2007}, with a
mean measure $\nu$ (as discussed in detail in Section \ref{sec.BPasPP}). 
Therefore, the process presented in \cite{Paisley:2010} must also be a Poisson process 
with this same mean measure.  Showing this equivalence would provide a direct proof 
of \cite{Paisley:2010} using the well-studied Poisson process machinery \citep{Kingman:1993}.

In this paper we present such a derivation (Section~\ref{sec.BPproof}).  
In addition, we derive error truncation bounds that are 
tighter than those in the literature (Section~\ref{sec.truncations}) 
\citep{Teh:2009,Paisley:2011}.  The Poisson process framework also 
provides an immediate proof of the extension of the construction to 
beta processes with a varying concentration parameter and infinite base measure (Section~\ref{sec.extension}), which does not follow immediately from the derivation in \cite{Paisley:2010}. In Section \ref{sec.MCMC}, we present a new MCMC algorithm for stick-breaking beta processes that uses the Poisson process to yield a more efficient sampler than that presented in \cite{Paisley:2010}.

\begin{figure*}[t]
\centering
\includegraphics[width=.9\textwidth]{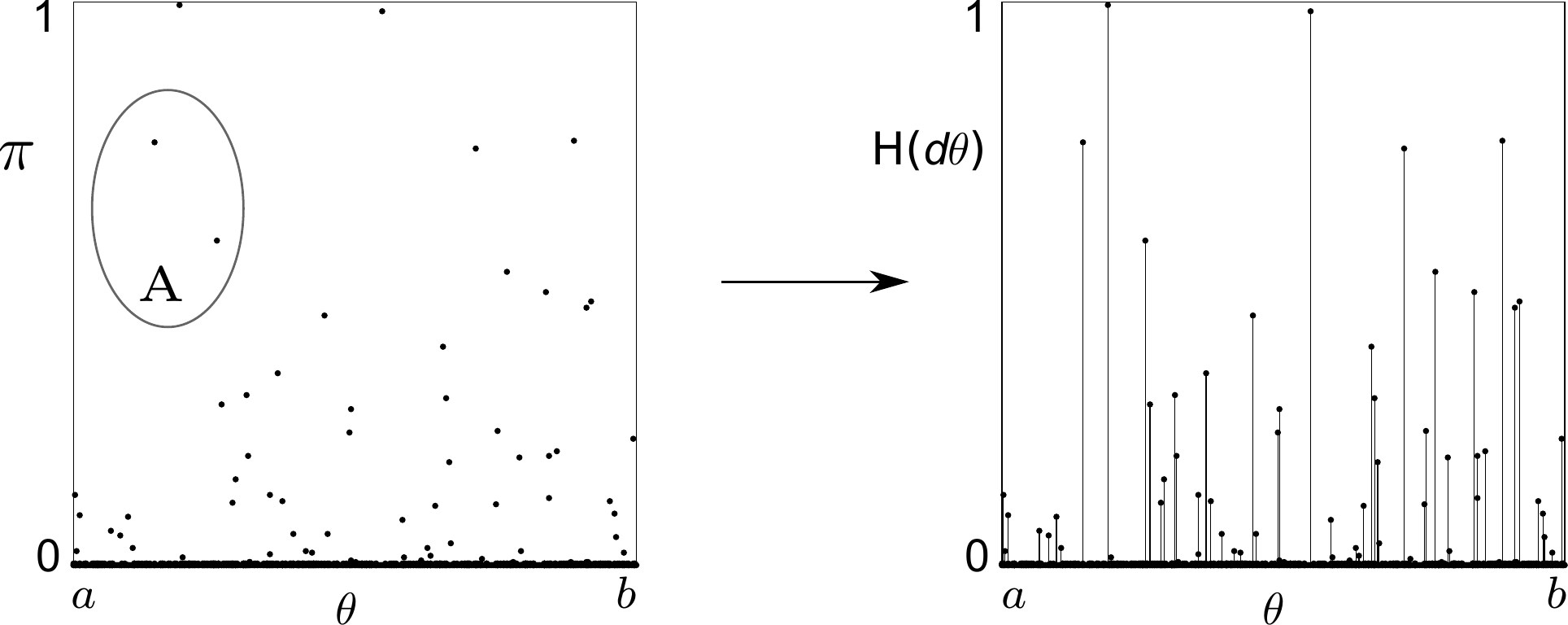}
\caption{(left) A Poisson process $\rmPi$ on $[a,b]\times[0,1]$ with mean measure $\nu = \mu \times \lambda$, where $\lambda(d\pi) = \alpha \pi^{-1}(1-\pi)^{\alpha-1}d\pi$ and $\mu([a,b]) < \infty$. The set $A$ contains a Poisson distributed number of atoms with parameter $\int_A \mu(d\theta)\lambda(d\pi)$. (right) The beta process constructed from $\rmPi$. The first dimension corresponds to location, and the second dimension to weight.}\label{fig.betaprocess}
\end{figure*}

\section{The Beta Process}
In this section, we review the beta process and its marginalized representation. We discuss the link between the beta process and the Poisson process, defining the underlying L\'{e}vy measure of the beta process. We then review the stick-breaking construction of the beta process, and give an equivalent representation of the generative process that will help us derive its L\'{e}vy measure.

A draw from a beta process is (with probability one) a countably infinite 
collection of weighted atoms in a space $\Omega$, with weights that lie in 
the interval $[0,1]$~\citep{Hjort:1990}. Two parameters govern the distribution on these weights, a concentration parameter $\alpha > 0$ and a finite base measure 
$\mu$, with $\mu(\Omega) = \gamma$.\footnote{In Section \ref{sec.extension} we discuss a generalization of this definition that is more in line with the definition given by \cite{Hjort:1990}.}
%\footnote{\cite{Hjort:1990} originally 
%developed the beta process in the context of survival analysis, where 
%$\Omega = \mathbb{R}_+$ and $\mu$ can be a $\sigma$-finite measure with 
%$\mu(\mathbb{R}_+) = \infty$, and $\alpha$ a time-varying function, $\alpha(\theta)$. 
%In Section \ref{sec.extension} we discuss how the construction in this paper 
%extends to this case.}  
Since such a draw is an atomic measure, we can write it as 
$H = \sum_{ij} \pi_{ij}\delta_{\theta_{ij}}$, where the two index values follow from \cite{Paisley:2010}, and we write 
$H \sim \mathrm{BP}(\alpha,\mu)$.

Contrary to the Dirichlet process, which provides a probability measure, the total measure $H(\Omega) \neq 1$ with probability one. Instead, beta processes are useful as parameters for a Bernoulli process. We write the Bernoulli process $X$ as $X = \sum_{ij}z_{ij}\delta_{\theta_{ij}}$, where 
$z_{ij} \sim \mathrm{Bernoulli}(\pi_{ij})$, and denote this as $X \sim \mathrm{BeP}(H)$.  \cite{Thibaux:2007} show that
marginalizing over $H$ yields the Indian buffet process (IBP) of \cite{Griffiths:2006}.  

The IBP clearly shows the featural clustering property of the beta process, and is specified as follows: To generate a sample $X_{n+1}$ from an IBP conditioned 
on the previous $n$ samples, draw \vspace{-1mm}
\begin{equation}
X_{n+1}|X_{1:n} \sim \mathrm{BeP}\left(\frac{1}{\alpha + n}\sum_{m=1}^n X_m  + \frac{\alpha}{\alpha + n} \mu\right).\nonumber
\end{equation}
This says that, for each $\theta_{ij}$ with at least one value of $X_m(\theta_{ij})$ equal
to one, the value of $X_{n+1}(\theta_{ij})$ is equal to one with probability $\frac{1}{\alpha + n}\sum_m X_m(\theta_{ij})$. After sampling these locations, a $\mathrm{Poisson}(\alpha\mu(\Omega)/(\alpha + n))$ distributed number of new locations $\theta_{i'j'}$ are introduced with corresponding $X_{n+1}(\theta_{i'j'})$ set equal to one. From this representation one can show that $X_m(\Omega)$ has a $\mathrm{Poisson}(\mu(\Omega))$ distribution, and the number of unique observed atoms in the process $X_{1:n}$ is Poisson distributed with parameter $\sum_{m=1}^n \alpha\mu(\Omega)/(\alpha + m - 1)$ \citep{Thibaux:2007}.

\subsection{The beta process as a Poisson process}\label{sec.BPasPP}
An informative perspective of the beta process is as a {completely
random measure}, a construction based on the Poisson 
process \citep{Jordan:2010}. We illustrate 
this in Figure \ref{fig.betaprocess} using an example where 
$\Omega = [a,b]$ and $\mu(A) = \frac{\gamma}{b-a} \mathrm{Leb}(A)$, with $\mathrm{Leb}(\cdot)$ the Lebesgue measure. 
The right figure shows a draw from the beta process. The left figure shows 
the underlying Poisson process, $\rmPi = \{(\theta,\pi)\}$. 

In this example, a Poisson process generates points in the space $[a,b]\times[0,1]$. 
It is completely characterized by its {mean measure},
$\nu(d\theta,d\pi)$~\citep{Kingman:1993,Cinlar:2011}. For any subset $A\subset [a,b]\times[0,1]$, the random {counting measure} $N(A)$ equals the number 
of points from $\rmPi$ contained in $A$. The distribution of 
$N(A)$ is Poisson with parameter $\nu(A)$.  Moreover, for all pairwise disjoint sets 
$A_1,\dots,A_n$, the random variables $N(A_1),\dots,N(A_n)$ are independent, 
and therefore $N$ is completely random. 

In the case of the beta process, the mean measure of the underlying Poisson process is
\begin{equation}\label{eqn.PPmeanmeasure}
\nu(d\theta,d\pi) = \alpha \pi^{-1}(1-\pi)^{\alpha-1} d\pi \mu(d\theta).
\end{equation}
We refer to $\lambda(d\pi) = \alpha \pi^{-1}(1-\pi)^{\alpha-1} d\pi$ as the 
{L\'{e}vy measure} of the process, and $\mu$ as its base measure. Our goal in Section \ref{sec.BPproof} 
will be to show that the following construction is also a Poisson process with mean measure equal 
to (\ref{eqn.PPmeanmeasure}), and is thus a beta process.

\begin{figure*}[t]
 \centering
\includegraphics[width=.9\textwidth]{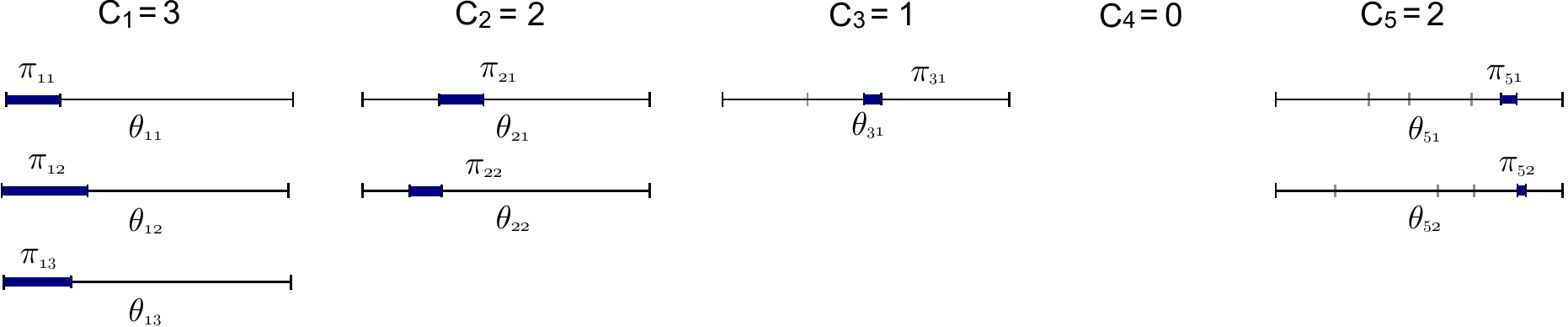}
\caption{An illustration of the stick-breaking construction of the beta process by round index $i$ for $i\leq 5$. Given a space $\Omega$ with measure $\mu$, for each index $i$ a $\mathrm{Poisson}(\mu(\Omega))$ distributed number of atoms $\theta$ are drawn i.i.d.\ from the probability measure $\mu/\mu(\Omega)$. To atom $\theta_{ij}$, a corresponding weight $\pi_{ij}$ is attached that is the $i$th break drawn independently from a $\mathrm{Beta}(1,\alpha)$ stick-breaking process. A beta process is $H = \sum_{ij}\pi_{ij}\delta_{{\theta}_{ij}}$.}\label{fig.stickpic}
\end{figure*}

\subsection{Stick-breaking for the beta process}\label{sec.stickbreaking}
\cite{Paisley:2010} presented a method for explicitly constructing beta processes
based on the notion of stick-breaking, a general method for obtaining discrete 
probability measures \citep{Ishwaran:2001}.  Stick-breaking plays an important
role in Bayesian nonparametrics, thanks largely to a seminal derivation of a stick-breaking
representation for the Dirichlet process by \cite{Sethuraman:1994}.  In the case
of the beta process, \cite{Paisley:2010} presented the following representation:
\begin{equation}\label{eqn.Horiginal}
 H = \sum_{i=1}^{\infty} \sum_{j=1}^{C_i} V_{ij}^{(i)}\prod_{l=1}^{i-1}(1-V_{ij}^{(l)})\delta_{\theta_{ij}},
\end{equation}
\begin{equation}\nonumber
C_i \iid \mathrm{Poisson}(\gamma),\quad V_{ij}^{(l)} \iid \mathrm{Beta}(1,\alpha),\quad \theta_{ij} \iid \frac{1}{\gamma}\mu ,
\end{equation}
where, as previously mentioned, $\alpha > 0$ and $\mu$ is a non-atomic finite base measure with 
$\mu(\Omega) = \gamma$.  

This construction sequentially incorporates into $H$ a Poisson-distributed number of atoms drawn i.i.d.\ from $\mu/\gamma$, with each round in this sequence indexed by $i$. The atoms receive weights in $[0,1]$, drawn independently according to a stick-breaking construction---an atom in round $i$ throws away the first $i-1$ breaks from its stick, and keeps the $i$th break as its weight. We illustrate this in Figure \ref{fig.stickpic}.

We use an equivalent definition of $H$ that reduces the total number of random variables by reducing the product $\prod_{j<i}(1-V_j)$ to a function of a single random variable. Let $V_i$ be i.i.d.\ $\mathrm{Beta}(1,\alpha)$ and let $f(V_{1:i-1}) := \prod_{j<i}(1-V_j)$. If $T \sim \mathrm{Gamma}(i-1,\alpha)$, then $f(V_{1:i-1}) =_d \exp\{-T\}$.
%\footnote{We derive this as follows: After a change of variables, $-\ln (1-V_j) \sim \mathrm{Exponential}(\alpha)$. The sum of $i-1$ i.i.d.\ exponential random variables is $\mathrm{Gamma}(i-1,\alpha)$ distributed, and thus $T =_d \sum_{j=1}^{i-1}-\ln (1-V_j)$. The equivalence of (\ref{eqn.H}) follows.} 
The construction in (\ref{eqn.Horiginal}) is therefore equivalent to
\begin{equation}\label{eqn.H}\nonumber
 H = \sum_{j=1}^{C_1}V_{1j}\delta_{\theta_{1j}} + \sum_{i=2}^{\infty} \sum_{j=1}^{C_i} V_{ij}\mathrm{e}^{-T_{ij}}\delta_{\theta_{ij}},
\end{equation}
\begin{equation}\nonumber
C_i \iid \mathrm{Poisson}(\gamma),\quad V_{ij} \iid \mathrm{Beta}(1,\alpha),
\end{equation}
\begin{equation}
 T_{ij} \stackrel{ind}{\sim} \mathrm{Gamma}(i-1,\alpha),\quad \theta_{ij} \iid \frac{1}{\gamma}\mu .
\end{equation}
Starting from a finite approximation of the beta process, \cite{Paisley:2010}
showed that (\ref{eqn.Horiginal}) must be a beta process by making use of the stick-breaking construction of a beta distribution \citep{Sethuraman:1994}, and then finding the limiting case; a similar limiting-case derivation was given for the Indian buffet
process \citep{Griffiths:2006}.  We next show that (\ref{eqn.Horiginal})
can be derived directly from the characterization of the beta process as a
Poisson process. This verifies the construction, and also leads to new properties of the beta process.

\section{Stick-breaking from the Poisson Process}\label{sec.BPproof}
We now prove that (\ref{eqn.Horiginal}) is a beta process with parameter
$\alpha > 0$ and base measure $\mu$ by showing that it has an underlying Poisson process with mean measure (\ref{eqn.PPmeanmeasure}).\footnote{A similar
result has recently been presented by \cite{Broderick:2011}; however, their approach differs 
from ours in its mathematical underpinnings.  Specifically we use a decomposition of the beta 
process into a countably infinite collection of Poisson processes, which leads directly 
to the applications that we pursue in subsequent sections.
By contrast, the proof in \cite{Broderick:2011} does not take this route, and their focus is on power-law generalizations of the beta process.}
We first state two basic lemmas regarding Poisson processes \citep{Kingman:1993}.
We then use these lemmas to show that the construction of $H$ in
(\ref{eqn.H}) has an underlying Poisson process representation, followed by the proof.

\subsection{Representing $H$ as a Poisson process}
The first lemma concerns the marking of points in a Poisson process with i.i.d.\ random variables. The second lemma concerns the superposition of independent Poisson processes. Theorem 1 uses these two lemmas to show that the construction in (\ref{eqn.H}) has an underlying Poisson process.

\paragraph{Lemma 1 (marked Poisson process)} \emph{Let $\rmPi^*$ be a Poisson process on $\Omega$ with mean measure $\mu$. With each $\theta \in \rmPi^*$ associate a random variable $\pi$ drawn independently with probability measure $\lambda$ on $[0,1]$. Then the set $\rmPi = \{(\theta,\pi)\}$ is a Poisson process on $\Omega \times [0,1]$ with mean measure $\mu \times \lambda$.}

\paragraph{Lemma 2 (superposition property)} \emph{Let $\rmPi_1, \rmPi_2,\dots$ be a countable collection of independent Poisson processes on $\Omega\times [0,1]$. Let $\rmPi_i$ have mean measure $\nu_i$. Then the superposition $\rmPi = \bigcup_{i=1}^{\infty} \rmPi_i$ is a Poisson process with mean measure $\nu = \sum_{i=1}^{\infty} \nu_i$.}

\paragraph{Theorem 1} \emph{The construction of $H$ given in (\ref{eqn.H}) has an underlying Poisson process.}

\emph{Proof.} This is an application of Lemmas 1 and 2; in this proof we fix some notation for what follows. Let $\pi_{1j} := V_{1j}$ and $\pi_{ij} := V_{ij}\exp\{-T_{ij}\}$ for $i > 1$. Let $H_i := \sum_{j=1}^{C_i} \pi_{ij}\delta_{\theta_{ij}}$ and therefore $H = \sum_{i=1}^{\infty} H_i$. Noting that $C_i \sim \mathrm{Poisson}(\mu(\Omega))$, for each $H_i$ the set of atoms $\{\theta_{ij}\}$ forms a Poisson process $\Pi^*$ on $\Omega$ with mean measure $\mu$. Each $\theta_{ij}$ is marked with a $\pi_{ij} \in [0,1]$ that has some probability measure $\lambda_i$ (to be defined later).  By Lemma 1, each $H_i$ has an underlying Poisson process $\rmPi_i = \{(\theta_{ij},\pi_{ij})\}$, on $\Omega \times [0,1]$ with mean measure $\mu\times \lambda_i$. It follows that $H$ has an underlying $\rmPi = \bigcup_{i=1}^{\infty} \rmPi_i$, which is a superposition of a countable collection of independent Poisson processes, and is therefore a Poisson process by Lemma 2.$\hfill\square$

\subsection{Calculating the mean measure of $H$}\label{sec.proof}
We've shown that $H$ has an underlying Poisson process; it remains to calculate its mean measure. We define the mean measure of $\rmPi_i$ to be $\nu_i = \mu\times \lambda_i$, and by Lemma 2 the mean measure of $\rmPi$ is $\nu = \sum_{i=1}^{\infty} \nu_i = \mu\times\sum_{i=1}^{\infty} \lambda_i$. We next show that $\nu(d\theta,d\pi) = \alpha \pi^{-1}(1-\pi)^{\alpha-1}d\pi\mu(d\theta)$, which will establish the result stated in the following theorem.

\paragraph{Theorem 2} \emph{The construction defined in (\ref{eqn.Horiginal}) is of a beta process with parameter $\alpha > 0$ and finite base measure $\mu$.}

\noindent \emph{Proof.} To show that the mean measure of $\rmPi$ is equal to (\ref{eqn.PPmeanmeasure}), we first calculate each $\nu_i$ and then take their summation. We split this calculation into two groups, $\rmPi_1$ and $\rmPi_i$ for $i > 1$, since the distribution of $\pi_{ij}$ (as defined in the proof of Theorem 1) requires different calculations for these two groups. We use the definition of $H$ in (\ref{eqn.H}) to calculate these distributions of $\pi_{ij}$ for $i > 1$.

\emph{Case} $i=1$.~~ The first round of atoms and their corresponding weights, $H_1 = \sum_{j=1}^{C_1} \pi_{1j}\delta_{\theta_{1j}}$ with $\pi_{1j} := V_{1j}$, has an underlying Poisson process $\rmPi_1 = \{(\theta_{1j},\pi_{1j})\}$ with mean measure $\nu_1 = \mu \times \lambda_1$ (Lemma 1). It follows that
\begin{equation}\label{eqn.f_1}
 \lambda_1(d\pi) = \alpha (1-\pi)^{\alpha-1}d\pi.
\end{equation}
We write $\lambda_i(d\pi) = f_i(\pi|\alpha)d\pi$. For example, the density above is $f_1 = \alpha(1-\pi)^{\alpha-1}$. We next focus on calculating the density $f_i$ for $i>1$.

\emph{Case} $i > 1$.~~ Each $H_i$ has an underlying Poisson process $\rmPi_i = \{(\theta_{ij},\pi_{ij})\}$ with mean measure $\mu\times\lambda_i$, where $\lambda_i$ determines the probability distribution of $\pi_{ij}$ (Lemma 1). As with $i=1$, we write this measure as $\lambda_i(d\pi) = f_i(\pi|\alpha)d\pi$, where $f_i(\pi|\alpha)$ is the density of $\pi_{ij}$, i.e., of the $i$th break from a $\mathrm{Beta}(1,\alpha)$ stick-breaking process. This density plays a significant role in the truncation bounds and MCMC sampler derived in the following sections; we next focus on its derivation.

Recall that $\pi_{ij} := V_{ij}\exp\{-T_{ij}\}$, where $V_{ij} \sim \mathrm{Beta}(1,\alpha)$ and $T_{ij} \sim \mathrm{Gamma}(i-1,\alpha)$. First, let $W_{ij} := \exp\{-T_{ij}\}$. Then by a change of variables,
$$p_W(w|i,\alpha) = \frac{\alpha^{i-1}}{(i-2)!} w^{\alpha-1}(-\ln w)^{i-2}\,.$$
Using the product distribution formula for two random variables \citep{Rohatgi:1976}, the density of $\pi_{ij} = V_{ij}W_{ij}$ is
\begin{eqnarray}\label{eqn.f_i}
f_i(\pi|\alpha) \hspace{-2mm}&=&\hspace{-2mm} \int_{\pi}^1 w^{-1}p_V(\pi/w|\alpha)p_W(w|i,\alpha)dw\\
\hspace{-2mm}&=&\hspace{-2mm} \frac{\alpha^{i}}{(i-2)!}\int_{\pi}^1 w^{\alpha-2}(\ln \frac{1}{w})^{i-2}(1-\frac{\pi}{w})^{\alpha-1} dw\nonumber.
\end{eqnarray}
Though this integral does not have a closed-form solution for a single L\'{e}vy measure $\lambda_i$, we show next that the sum over these measures does have a closed-form solution. 

\emph{The L\'{e}vy measure of $H$.}~~ Using the values of $f_i$ derived above, we can calculate the mean measure of the Poisson process underlying (\ref{eqn.Horiginal}). As discussed, the measure $\nu$ can be decomposed as follows,
$$\nu(d\theta,d\pi) = \sum_{i=1}^{\infty} (\mu\times \lambda_i)(d\theta, d\pi) = \mu(d\theta)d\pi\sum_{i=1}^{\infty} f_i(\pi|\alpha).$$
By showing that $\sum_{i=1}^{\infty} f_i(\pi|\alpha) = \alpha\pi^{-1}(1-\pi)^{\alpha-1}$, we complete the proof; we refer to the appendix for the details of this calculation.

\section{Some Properties of the Beta Process}\label{sec.discussion}
We have shown that the stick-breaking construction defined in 
(\ref{eqn.Horiginal}) has an underlying Poisson process with mean measure 
$\nu(d\theta,d\pi) = \alpha \pi^{-1}(1-\pi)^{\alpha-1}d\pi\mu(d\theta)$, 
and is therefore a beta process.  Representing the stick-breaking 
construction as a superposition of a countably infinite collection of 
independent Poisson processes is also useful for further characterizing 
the beta process. For example, we can use this representation to analyze 
truncation properties. We can also easily extend the construction in 
(\ref{eqn.Horiginal}) to cases such as that considered in \cite{Hjort:1990}, 
where $\alpha$ is a function of $\theta$ and $\mu$ is an infinite measure.

\begin{figure*}
\centering
\includegraphics[width=.9\textwidth]{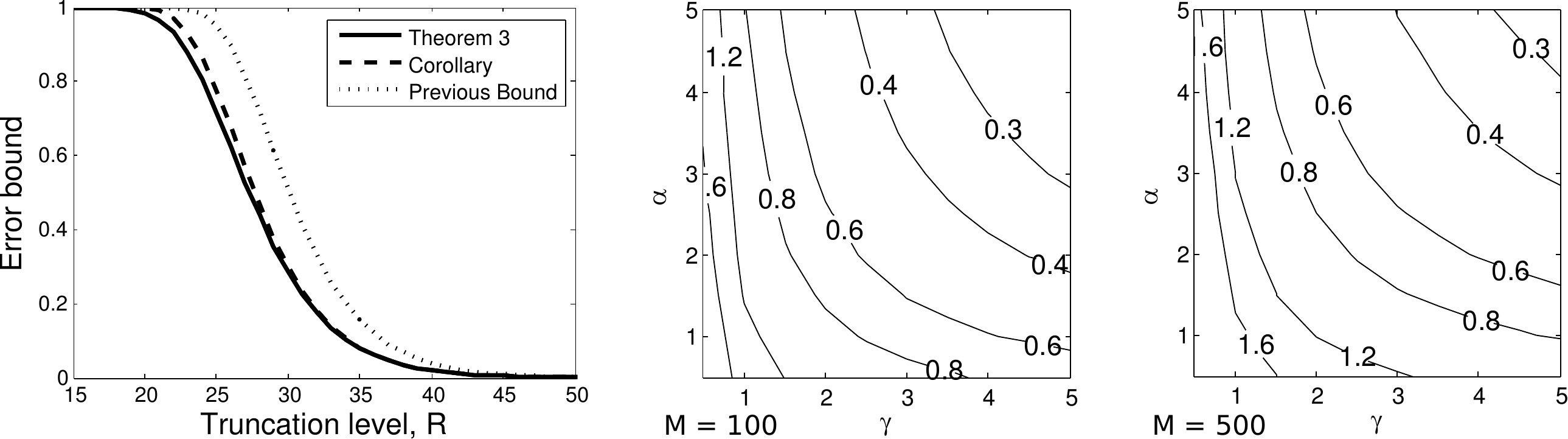}
\caption{Examples of the error bound. (left) The bounds for $\alpha = 3$, $\gamma = 4$ and $M = 500$. The previous bound appears in \cite{Paisley:2011}. (center and right) Contour plots of the $L_1$ distance between the Theorem 3 bound and the Corollary bound, presented as functions of $\alpha$ and $\gamma$ for (center) $M=100$, (right) $M=500$. The $L_1$ distance for the left plot is $0.46$. The Corollary bound becomes tighter as $\alpha$ and $\gamma$ increase, and as $M$ decreases.}\label{fig.bound}
\end{figure*}

\subsection{Truncated beta processes}\label{sec.truncations}
Truncated beta processes arise in the variational inference setting 
\citep{Teh:2009,Paisley:2011,Jordan:1999}. Poisson process 
representations are useful for characterizing the part of the beta process 
that is being thrown away in the truncation. Consider a beta process 
truncated after round $R$, defined as $H^{(R)} = \sum_{i=1}^R H_i$. 
The part being discarded, $H - H^{(R)}$, has an underlying Poisson process 
with mean measure 
\begin{eqnarray}
\nu_R^+(d\theta,d\pi) & := & \textstyle\sum_{i=R+1}^{\infty} \nu_i(d\theta,d\pi)\nonumber\\
		      & = & \mu(d\theta)\times\textstyle\sum_{i=R+1}^{\infty} \lambda_i(d\pi),
\end{eqnarray}
and a corresponding counting measure $N_R^+(d\theta,d\pi)$. This measure 
contains information about the missing atoms.\footnote{For example, the number of 
missing atoms having weight $\pi \geq \epsilon$ is Poisson distributed with 
parameter $\nu_R^+(\Omega,[\epsilon,1])$.}

% while the expected number of 
% missing ones contained in the discarded portion of a Bernoulli process, 
% $\mathbb{E}[X_R^+(\Omega)]$ where $X_R^+\sim\mathrm{BeP}(H-H^{(R)})$, is equal 
% to $\mathbb{E}\int_0^1 \pi N_R^+(\Omega,d\pi) = \int_0^1 \pi \nu_R^+(\Omega,d\pi)$.

For truncated beta processes, a measure of closeness to the true beta process is helpful when selecting truncation levels. To this end, let data $Y_n \sim f(X_n,\phi_n)$, where $X_n$ is a Bernoulli process taking either $H$ or $H^{(R)}$ as parameters, and $\phi_n$ is a set of additional parameters (which could be globally shared). Let $\Y = (Y_1,\dots,Y_M)$. One measure of closeness is the $L_1$ distance between the marginal density of $\Y$ under the beta process, $\m_{\infty}(\Y)$, and the process truncated at round $R$, $\m_R(\Y)$. This measure originated with work on truncated Dirichlet processes in \cite{Ishwaran:2000,Ishwaran:2001}; in \cite{Teh:2009}, it was extended to the beta process. 

After slight modification to account for truncating rounds rather than atoms, 
the result in \cite{Teh:2009} implies that
\begin{equation}\label{eqn.bound}
\frac{1}{4}\int |\m_R(\Y) - \m_{\infty}(\Y)|d\Y  \hspace{3cm}
\end{equation}
$$\quad\quad\quad \leq  \mathbb{P}\left(\exists (i,j), i > R, 1 \leq n \leq M : X_n(\theta_{ij}) = 1\right),$$
with a similar proof as in \cite{Ishwaran:2000}. This says that 1/4 times the 
$L_1$ distance between $\m_R$ and $\m_{\infty}$ is less than one minus the 
probability that, in $M$ Bernoulli processes with parameter 
$H\sim\mathrm{BP}(\alpha,\mu)$, there is no $X_n(\theta) = 1$ for a $\theta \in H_i$ with $i> R$. 
In \cite{Teh:2009} and \cite{Paisley:2011}, this bound was loosened. Using the Poisson process 
representation of $H$, we can give an exact form of this bound. To do so, we use 
the following lemma, which is similar to Lemma 1, but accounts for markings 
that are not independent of the atom.

\noindent\textbf{Lemma 3} \emph{Let $(\theta,\pi)$ form a Poisson process on $\Omega\times [0,1]$ with mean measure $\nu_R^+$. Mark each $(\theta,\pi)$ with a random variable $U$ in a finite space $\mathcal{S}$ with transition probability kernel $Q(\pi,\cdot)$. Then $(\theta,\pi,U)$ forms a Poisson process on $\Omega\times [0,1]\times \mathcal{S}$ with mean measure $\nu_R^+(d\theta,d\pi)Q(\pi,U)$.}

This leads to Theorem 3.

\paragraph{Theorem 3} \emph{Let $X_{1:M} \iid \mathrm{BeP}(H)$ with $H \sim \mathrm{BP}(\alpha,\mu)$ constructed as in (\ref{eqn.Horiginal}). For a truncation value $R$, let $E$ be the event that there exists an index $(i,j)$ with $i>R$ such that $X_n(\theta_{ij}) = 1$. Then the bound in (\ref{eqn.bound}) equals}
$$\mathbb{P}(E) = 1 - \exp\left\lbrace - \int_{(0,1]} \nu_R^+(\Omega,d\pi)\left(1 - (1-\pi)^M\right)\right\rbrace.$$
\emph{Proof.} Let $U\in\{0,1\}^M$. By Lemma 3, the set $\{(\theta,\pi,U)\}$ constructed from rounds $R+1$ and higher is a Poisson process on $\Omega\times [0,1]\times\{0,1\}^M$ with mean measure $\nu_R^+(d\theta,d\pi)Q(\pi,U)$ and a corresponding counting measure $N_R^+(d\theta,d\pi,U)$, where $Q(\pi,\cdot)$ is a transition probability measure on the space $\{0,1\}^M$. Let $A = \{0,1\}^M \backslash \boldsymbol{0}$, where $\boldsymbol{0}$ is the zero vector. Then $Q(\pi,A)$ is the probability of this set with respect to a Bernoulli process with parameter $\pi$, and therefore $Q(\pi,A) = 1 - (1-\pi)^M$. The probability $\mathbb{P}(E) = 1 - \mathbb{P}(E^c)$, which is equal to $1 - \mathbb{P}(N_R^+(\Omega,[0,1],A)=0)$. The theorem follows since $N_R^+(\Omega,[0,1],A)$ is a Poisson-distributed random variable with parameter $\int_{(0,1]}\nu_R^+(\Omega,d\pi)Q(\pi,A)$.\footnote{We give a second proof using simple functions in the appendix. One can use approximating simple functions to give an arbitrarily close approximation of Theorem 3. Furthermore, since $\nu_R^+ = \nu_{R-1}^+ - \nu_{R}$ and $\nu_0^+ = \nu$, performing a sweep of truncation values requires approximating only one additional integral for each increment of $R$. }$\hfill\square$ 

%By exchanging integrals with respect to $\pi$ and $w$, only a one-dimensional integral needs approximating.

Using the Poisson process, we can give an analytical bound that is tighter than that in \cite{Paisley:2011}.
\paragraph{Corollary 1} \emph{Given the set-up in Theorem 3, an upper bound on $\mathbb{P}(E)$ is} $$\mathbb{P}(E) \leq 1 - \exp\left\lbrace-\gamma M \left(\frac{\alpha}{1+\alpha}\right)^R\right\rbrace.$$

\emph{Proof.} We give the proof in the appendix.

The bound in \cite{Paisley:2011} has $2M$ rather than $M$. We observe that the term in the exponential equals the negative of $M\int_0^1 \pi \nu_R^+(\Omega,d\pi)$, which is the expected number of missing ones in $M$ truncated Bernoulli process observations. Figure \ref{fig.bound} shows an example of these bounds. %Note that, for a value of $R$, the number of observed atoms has a $\mathrm{Poisson}(\gamma R)$ distribution.

\subsection{Beta processes with infinite $\mu$ and varying $\alpha$}\label{sec.extension}
The Poisson process allows for the construction to be extended to the more general definition of the beta process given by \cite{Hjort:1990}.
In this definition, the value of 
$\alpha(\theta)$ is a function of $\theta$, rather than a constant, and the base measure $\mu$ may be infinite, but $\sigma$-finite.\footnote{That is, the total measure
$\mu(\Omega) = \infty$, but there is a measurable partition $(E_k)$ of $\Omega$ with each $\mu(E_k) < \infty$.}
Using Poisson processes, the extension of (\ref{eqn.Horiginal}) to this setting is straightforward. We note that this is not immediate from the limiting case derivation presented in \cite{Paisley:2010}.

For a partition $(E_k)$ of $\Omega$ with $\mu(E_k) < \infty$, we treat each set $E_k$ as a separate Poisson 
process with mean measure
\begin{eqnarray}\nonumber
\nu_{E_k}(d\theta,d\pi) &=& \mu(d\theta)\lambda(\theta,d\pi),\quad \theta\in E_k\\
&=& \alpha(\theta)\pi^{-1}(1-\pi)^{\alpha(\theta)-1}d\pi\mu(d\theta).\nonumber
\end{eqnarray}
The transition probability 
kernel $\lambda$ follows from the continuous version of Lemma 3. By superposition, we have the overall beta
process.  Modifying (\ref{eqn.Horiginal}) gives the following construction: For each set $E_k$ construct a separate $H_{E_k}$.
In each round of (\ref{eqn.Horiginal}), incorporate $\mathrm{Poisson}(\mu(E_k))$ 
new atoms $\theta_{ij}^{(k)}\in E_k$ drawn i.i.d.\ from $\mu/\mu(E_k)$. For atom $\theta_{ij}^{(k)}$, draw a weight $\pi_{ij}^{(k)}$ using the $i$th break 
from a $\mathrm{Beta}(1,\alpha(\theta_{ij}^{(k)}))$ stick-breaking process. 
The complete beta process is the union of these local beta processes.

\section{MCMC Inference}\label{sec.MCMC}
We derive a new MCMC inference algorithm for beta processes that incorporates ideas from the stick-breaking construction and Poisson process. In the algorithm, we re-index atoms to take one index value $k$, and let $d_k$ indicate the Poisson process of the $k$th atom under consideration (i.e., $\theta_{k} \in H_{d_k}$). For calculation of the likelihood, given $M$ Bernoulli process draws, we denote the sufficient statistics $m_{1,k} = \sum_{n=1}^M X_n(\theta_k)$ and $m_{0,k} = M - m_{1,k}$.

We use the densities $f_1$ and $f_i$, $i > 1$, derived in (\ref{eqn.f_1}) and (\ref{eqn.f_i}) above. Since the numerical integration in (\ref{eqn.f_i}) is computationally expensive, we sample $w$ as an auxiliary variable. The joint density of $\pi_{ij}$ and $w_{ij}$, $0 \leq \pi_{ij} 
\leq w_{ij}$, for $\theta_{ij} \in H_i$ and $i>1$ is
\begin{equation}\nonumber
f_i(\pi_{ij},w_{ij}|\alpha) \propto w_{ij}^{-1}(-\ln w_{ij})^{i-2}(w_{ij}-\pi_{ij})^{\alpha-1}.
\end{equation}
The density for $i=1$ does not depend on $w$.

\subsection{A distribution on observed atoms}\label{sec.observed}
Before presenting the MCMC sampler, we derive a quantity that we use in the algorithm. Specifically, for the collection of Poisson processes $H_i$, we calculate the distribution on the number of atoms $\theta \in H_i$ for which the Bernoulli process $X_n(\theta)$ is equal to one for some $1 \leq n \leq M$. In this case, we denote the atom as being ``observed.'' This distribution is relevant to inference, since in practice we care most about samples at these locations.

The distribution of this quantity is related to Theorem 3. There, the exponential term gives the probability that this number is zero for all $i > R$. More generally, under the prior on a single $H_i$, the number of observed atoms is Poisson distributed with parameter
\begin{equation}\label{eqn.nonzero}
\xi_i = \int_0^1 \nu_i(\Omega,d\pi)(1-(1-\pi)^M)d\pi.
\end{equation}
The sum $\sum_{i=1}^{\infty} \xi_i < \infty$ for finite $M$, meaning a finite number of atoms will be observed with probability one.

Conditioning on there being $T$ observed atoms overall, $\theta^*_{1:T}$, we can calculate a distribution on the Poisson process to which atom $\theta^*_k$ belongs. This is an instance of Poissonization of the multinomial; since for each $H_i$ the distribution on the number of observed atoms is independent and $\mathrm{Poisson}(\xi_i)$ distributed, conditioning on $T$ the Poisson process to which atom $\theta^*_k$ belongs is independent of all other atoms, and identically distributed with $\mathbb{P}(\theta^*_k \in H_i) \propto \xi_i$.

\subsection{The sampling algorithm}
We next present the MCMC sampling algorithm. We index samples by an $s$, and define all densities to be zero outside of their support.

\paragraph{Sample $\pi_k.$} We take several random walk Metropolis-Hastings steps for $\pi_k$. Let $\pi_k^s$ be the value at step $s$. 
Let the proposal be $\pi_k^{\star} = \pi_k^s + \xi_k^s$, where $\xi_k^s \iid N(0,\sigma^2_{\pi})$. Set 
$\pi_k^{s+1} = \pi_k^{\star}$ with probability $$\min \left\lbrace 1,\frac{p(m_{1,k},m_{0,k}|\pi_k^{\star})f_{d_k^s}(\pi_k^{\star}|w_k^s,\alpha_s)}{p(m_{1,k},m_{0,k}|\pi_k^s)f_{d_k^s}(\pi_k^s|w_k^s,\alpha_s)}\right\rbrace ,$$
otherwise set $\pi_k^{s+1} = \pi_k^{s}$. The likelihood and priors are
\begin{eqnarray}\nonumber
p(m_{1,k},m_{0,k}|\pi) &=& \pi^{m_{1,k}} (1-\pi)^{m_{0,k}},\\\nonumber
f_{d_k^s}(\pi|w_k^s,\alpha_s) &\propto& \left\{
\begin{array}{l l}
(w_k^s-\pi)^{\alpha_s-1}& \quad \mbox{if $d_k > 1$}\\
(1-\pi)^{\alpha_s-1}& \quad \mbox{if $d_k = 1$}
\end{array}\right .  .
\end{eqnarray}

\paragraph{Sample $w_k.$} We take several random walk Metropolis-Hastings steps for $w_k$ when $d_k > 1$. Let $w_k^s$ be the value at step $s$. Set the proposal $w_k^{\star} = w_k^s + \zeta_k^s$, where $\zeta_k^s \iid N(0,\sigma^2_w)$, and set $$w_k^{s+1} = w_k^{\star}\quad \mbox{w.p.}\quad \min \left\lbrace 1,\frac{f(w_k^{\star}|\pi_k^s,d_k^s,\alpha_s)}{f(w_k^s|\pi_k^s,d_k^s,\alpha_s)}\right\rbrace , $$
otherwise set $w_k^{s+1} = w_k^s$. The value of $f$ is
\begin{equation}\nonumber
f(w|\pi_k^s,d_k^s,\alpha_s) = w^{-1}(-\ln w)^{d_k^s-2}(w-\pi_k^s)^{\alpha_s-1} .
\end{equation}
When $d_k^s = 1$, the auxiliary variable $w_k$ does not exist, so we don't sample it. If $d_k^{s-1} = 1$, but $d_k^s > 1$, we sample $w_k^s \sim \mbox{Uniform}(\pi_k^s,1)$ and take many random walk M-H steps as detailed above.

\paragraph{Sample $d_k.$} We follow the discussion in Section \ref{sec.observed} to sample $d_k^{s+1}$. Conditioned on there being $T_s$ observed atoms at step $s$, the prior on $d_k^{s+1}$ is independent of all other indicators $d$, and $\mathbb{P}(d_k^{s+1} = i) \propto \xi_i^s$, where $\xi^s_i$ is given in (\ref{eqn.nonzero}). The likelihood depends on the current value of $d_k^s$.

\emph{Case} $d_k^s > 1.$ ~~The likelihood $f(\pi_k^s,w_k^s|d_k^{s+1}=i,\alpha_s)$ is proportional to
\begin{equation}\nonumber
\left\{
\begin{array}{l l}
\frac{\alpha_s^i}{(i-2)!} (w_k^s)^{-1}(-\ln w_k^s)^{i-2}(w_k^s-\pi_k^s)^{\alpha_s-1}&  \mbox{if $i > 1$}\vspace{2mm}\\
\alpha(1-\pi_k^s)^{\alpha_s-1}&  \mbox{if $i = 1$}
\end{array}\right .  
\end{equation}
\emph{Case} $d_k^s = 1.$ ~~In this case we must account for the possibility that $\pi_k^s$ may be greater than the most recent value of $w_k$, we marginalize the auxiliary variable $w$ numerically, and compute the likelihood as follows:
\begin{equation}\nonumber
\left\{
\begin{array}{l l}
\frac{\alpha_s^i}{(i-2)!}\textstyle\int_{\pi_k^s}^1 w^{-1}(-\ln w)^{i-2}(w-\pi_k^s)^{\alpha_s-1} dw&  \mbox{if $i > 1$}\vspace{2mm}\\
\alpha(1-\pi_k^s)^{\alpha_s-1}&  \mbox{if $i = 1$}
\end{array}\right .  
\end{equation}
A slice sampler \citep{Neal:2003} can be used to sample from this infinite-dimensional discrete distribution.

\paragraph{Sample $\alpha.$} We have the option of Gibbs sampling $\alpha$. For a $\mbox{Gamma}(\tau_1,\tau_2)$ prior on $\alpha$,  the full conditional of $\alpha$ is a gamma distribution with parameters
$$\tau_{1,s}' = \tau_1 + \textstyle\sum_k d_k^s,\quad\quad \tau_{2,s}' = \tau_2 - \textstyle\sum_k \ln (w_k^s - \pi_k^s).$$
In this case we set $w_k^s = 1$ if $d_k^s = 1$.

\paragraph{Sample $\gamma.$} We also have the option of Gibbs sampling $\gamma$ using a $\mbox{Gamma}(\kappa_1,\kappa_2)$ prior on $\gamma$. As discussed in Section \ref{sec.observed}, let $T_s$ be the number of observed atoms in the model at step $s$. The full conditional of $\gamma$ is a gamma distribution with parameters
$$\kappa_{1,s}' = \kappa_1 + T_s,\quad\quad \kappa_{2,s}' = \kappa_2 + \textstyle \sum_{n=0}^{M-1} \frac{\alpha_s}{\alpha_s + n}.$$
This distribution results from the Poisson process, and the fact that the observed and unobserved atoms form a disjoint set, and therefore can be treated as independent Poisson processes. In deriving this update, we use the equality $\sum_{i=1}^{\infty}\xi_i^s/\gamma_s = \sum_{n=0}^{M-1} \frac{\alpha_s}{\alpha_s + n}$, found by inserting the mean measure (\ref{eqn.PPmeanmeasure}) into (\ref{eqn.nonzero}).

\paragraph{Sample $X.$} For sampling the Bernoulli process $X$, we have that $p(X|\mathcal{D},H) \propto p(\mathcal{D}|X)p(X|H)$. The likelihood of data $\mathcal{D}$ is independent of $H$ given $X$ and is model-specific, while the prior on $X$ only depends on $\pi$.

\paragraph{Sample new atoms.} We sample new atoms in addition to the observed atoms. For each $i = 1,\dots,\mbox{max}(d_{1:T_s})$, we ``complete'' the round by sampling the unobserved atoms. For Poisson process $H_i$, this number has a $\mbox{Poisson}(\gamma_s - \xi_i^s)$ distribution. We can sample additional Poisson processes as well according to this distribution. In all cases, the new atoms are i.i.d.\ $\mu/\gamma_s$.

% \begin{figure*}[t]
% \centering
% \includegraphics[width=.99\textwidth]{BPimage.pdf}
% \caption{(left) A Poisson process $\rmPi$ on $[a,b]\times[0,1]$ with mean measure $\nu = \mu \times \lambda$, where $\lambda(d\pi) = \alpha \pi^{-1}(1-\pi)^{\alpha-1}d\pi$ and $\mu([a,b]) < \infty$. A set $A$ in this space contains a Poisson distributed number of atoms with parameter $\int_{(\theta,\pi)\in A} \mu(d\theta)\lambda(d\pi)$. (right) The beta process constructed from $\rmPi$. The first dimension corresponds to location, and the second dimension corresponds to weight.corresponds to weight.corresponds to weight.corresponds to weight.corresponds to weight.corresponds to weight.corresponds to weight.corresponds to weight.corresponds to weight.corresponds to weight.corresponds to weight.}
% \end{figure*}

% \begin{figure*}
% \centering
% \subfigure{\includegraphics[width=.25\textwidth]{toy_factors.pdf}}\quad\quad
% \subfigure{\includegraphics[width=.24\textwidth]{factor_hist.pdf}}\quad\quad
% \subfigure{\includegraphics[width=.29\textwidth]{density.pdf}}
% \caption{Results on synthetic data. (left) The top 16 underlying factor loadings for the MCMC iteration 10,000. The highly probable ``bars'' and ``diagonal'' patches are uncovered, as are the other randomly generated factors. (middle) A histogram of the number of factors. The empirical distribution centers on the truth. (right) The kernel smoothed density using the samples of $\alpha$ and $\gamma$. These values differ from the anticipated values (see the text for discussion).}\label{fig.results}
% \end{figure*}

\begin{figure*}
\centering
\includegraphics[width=.9\textwidth]{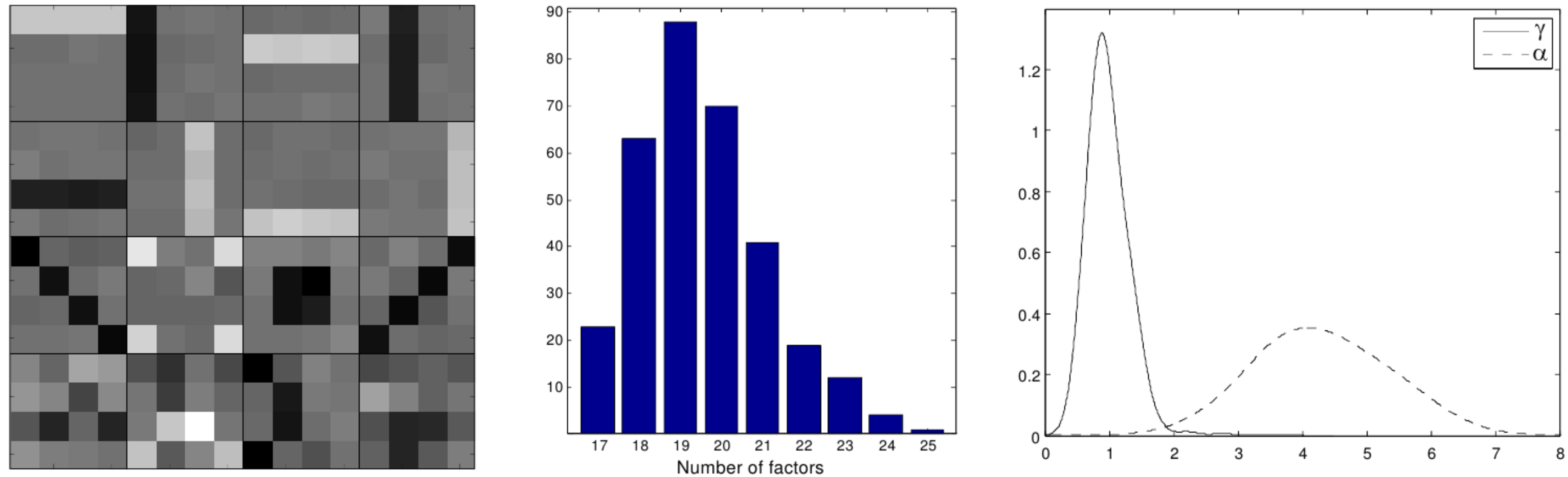}
\caption{Results on synthetic data. (left) The top 16 underlying factor loadings for MCMC iteration 10,000. The ground truth patterns are uncovered. (middle) A histogram of the number of factors. The empirical distribution centers on the truth. (right) The kernel smoothed density using the samples of $\alpha$ and $\gamma$ (see the text for discussion).}\label{fig.results}
\end{figure*}

\subsection{Experimental results}
We evaluate the MCMC sampler on synthetic data. We use the beta-Bernoulli process as a matrix factorization prior for a linear-Gaussian model. We generate a data matrix $Y = \Theta(W\circ Z) + \epsilon$ with each $W_{kn} \sim N(0,1)$, the binary matrix $Z$ has $\mbox{Pr}(Z_{kn}=1|H) = \pi_k$ and the columns of $\Theta$ are vectorized $4\times 4$ patches of various patterns (see Figure \ref{fig.results}). To generate $H$ for generating $Z$, we let $\pi_k$ be the expected value of the $k$th atom under the stick-breaking construction with parameters $\alpha = 1$, $\gamma = 2$. We place $\mbox{Gamma}(1,1)$ priors on $\alpha$ and $\gamma$ for inference. We sampled $M = 500$ observations, which used a total of 20 factors. Therefore $Y \in \mathrm{R}^{16\times 500}$ and $Z \in \{0,1\}^{20\times 500}$.

We ran our MCMC sampler for 10,000 iterations, collecting samples every 25th iteration after a burn-in of 2000 iterations. For sampling $\pi$ and $w$, we took 1,000 random walk steps using a Gaussian with variance $10^{-3}$. Inference was relatively fast; sampling all beta process related variables required roughly two seconds per iteration, which is significantly faster than the per-iteration average of 14 seconds for the algorithm presented in \cite{Paisley:2010}, where Monte Carlo integration was heavily used.

We show results in Figure \ref{fig.results}. While we expected to learn a $\gamma$ around two, and $\alpha$ around one, we note that our algorithm is inaccurate for these values. We believe that this is largely due to our prior on $d_k$ (Section \ref{sec.observed}). The value of $d_k$ significantly impacts the value of $\alpha$, and conditioning on $\sum_{n=1}^{M} X_n(\theta) > 0$ gives a prior for $d_k$ that is spread widely across the rounds and allows for much variation. A possible fix for this would be conditioning on the exact value of the number of atoms in a round. This will effectively give a unique prior for each atom, and would require significantly more numerical integrations leading to a slower algorithm.

Despite the inaccuracy in learning $\gamma$ and $\alpha$, the algorithm still found to the correct number of factors (initialized at 100), and found the correct underlying sparse structure of the data. This indicates that our MCMC sampler is able to perform the main task of finding a good sparse representation.\footnote{The variable $\gamma$ only enters the algorithm when sampling new atoms. Since we learn the correct number of factors, this indicates that our algorithm is not sensitive to $\gamma$. Fixing the concentration parameter $\alpha$ is an option, and is often done for Dirichlet processes.} It appeared that the likelihood of $\pi$ dominates inference for this value, since we observed that these samples tended to ``shadow'' the empirical distribution of $Z$.

\section{Conclusion}\label{sec.conclusion}
We have used the Poisson processes to prove that the stick-breaking construction presented by \cite{Paisley:2010} is a beta process. We then presented several consequences of this representation, including truncation bounds, a more general definition of the construction, and a new MCMC sampler for stick-breaking beta processes. Poisson processes offer flexible representations of Bayesian nonparametric priors; for example, \cite{Lin:2010} show how they can be used as a general representation of dependent Dirichlet processes. Representing a beta process as a superposition of a countable collection of Poisson processes may lead to similar generalizations.

\section*{Appendix}

\paragraph{Proof of Theorem 2 (conclusion)}
From the text, we have that $\lambda(d\pi) = f_1(\pi|\alpha)d\pi + \sum_{i=2}^{\infty}f_i(\pi|\alpha)d\pi$ with $f_1(\pi|\alpha) = \alpha(1-\pi)^{\alpha-1}$ and $f_i$ given 
in Equation \ref{eqn.f_i} for $i>1$. The sum of densities is
$$\textstyle\sum_{i=2}^{\infty} f_i(\pi|\alpha)\hspace{2.5in}\vspace{-.25cm}$$
\begin{eqnarray}\label{eqn.f_seq}
  &=& \sum_{i=2}^{\infty} \frac{\alpha^{i}}{(i-2)!}\int_{\pi}^1 w^{\alpha-2}(\ln \frac{1}{w})^{i-2}(1-\frac{\pi}{w})^{\alpha-1} dw\nonumber\\
&=&\hspace{-2mm}\alpha^2\int_{\pi}^1  w^{\alpha-2}(1-\frac{\pi}{w})^{\alpha-1}dw \sum_{i=2}^{\infty}\frac{\alpha^{i-2}}{(i-2)!}(\ln \frac{1}{w})^{i-2}\nonumber\\
&=&\hspace{-2mm}\alpha^2\int_{\pi}^1  w^{-2}(1-\pi/w)^{\alpha-1}dw\,.
\end{eqnarray}
The second equality is by monotone convergence and Fubini's theorem. This leads to an exponential power series, which simplifies to the third line. The last line equals $\frac{\alpha(1-\pi)^{\alpha}}{\pi}$. Adding the result of (\ref{eqn.f_seq}) to $\alpha(1-\pi)^{\alpha-1}$ gives $\sum_{i=1}^{\infty} f_i(\pi|\alpha) = \alpha\pi^{-1}(1-\pi)^{\alpha-1}.$ Therefore, $\nu(d\theta,d\pi) = \alpha\pi^{-1}(1-\pi)^{\alpha-1}d\pi\mu(d\theta)$, and the proof is complete. \hfill $\square$

\paragraph{Alternate proof of Theorem 3}
Let the set $B_{nk} = \left[\frac{k-1}{n},\frac{k}{n}\right)$ and $b_{nk} = \frac{k-1}{n}$, where $n$ and $k \leq n$ are positive integers.  Approximate the variable $\pi\in [0,1]$ with the simple function $g_n(\pi) = \sum_{k=1}^n b_{nk}\boldsymbol{1}_{B_{nk}}(\pi)$. We calculate the truncation error term, $\mathbb{P}(E^c)=\mathbb{E}[\prod_{i>R,j}(1-\pi_{ij})^M]$, by approximating with $g_n$, re-framing the problem as a Poisson process with mean and counting measures $\nu_R^+$ and $N_R^+(\Omega,B)$, and then taking a limit:

$\mathbb{E}\left[\prod\nolimits_{i>R,j}(1-\pi_{ij})^M\right] \hfill  $
\begin{eqnarray}
& = & \hspace{-1mm}\lim_{n\rightarrow\infty} \prod_{k=2}^n\mathbb{E}\left[(1-b_{nk})^{M\cdot N_R^+(\Omega,B_{nk})}\right]\\
%& = & \lim_{n\rightarrow\infty} \prod_{i=2}^n \exp\{-\nu_R^+(\Omega,B_{ni})\}\exp\{\nu_R^+(\Omega,B_{ni})(1-b_{ni})^M\}\nonumber\\
& = & \hspace{-1mm}\exp\left\lbrace \lim_{n\rightarrow\infty} - \sum_{k=2}^n \nu_R^+(\Omega,B_{nk})\left(1-(1-b_{nk})^M\right)\right\rbrace .\nonumber 
\end{eqnarray}
For a fixed $n$, this approach divides the interval $[0,1]$ into disjoint regions that can be analyzed separately as independent Poisson processes. Each region uses the approximation $\pi \approx g_n(\pi)$, with  $\lim_{n\rightarrow\infty}g_n(\pi) = \pi$, and $N_R^+(\Omega,B)$ counts the number of atoms with weights that fall in the interval $B$. Since $N_R^+$ is Poisson distributed with mean $\nu_R^+$, the expectation follows.

\paragraph{Proof of Corollary 1}
From the alternate proof of Theorem 3 above, we have $\mathbb{P}(E) = 1 - \mathbb{E}[\prod_{i>R,j}(1-\pi_{ij})^M] \leq 1 - \mathbb{E}[\prod_{i>R,j}(1-\pi_{ij})]^M$. This second expectation can be calculated as in Theorem 3 with $M$ replaced by a one. The resulting integral is analytic. Let $q_r$ be the distribution of the $r$th break from a $\mathrm{Beta}(1,\alpha)$ stick-breaking process. The negative of the term in the exponential of Theorem 3 is 
\begin{equation}\label{eqn.pfC1}
\int_0^1 \pi\nu_R^+(\Omega,d\pi) = \gamma \sum_{r=R+1}^{\infty}\mathbb{E}_{q_r}[\pi].
\end{equation}
Since $\mathbb{E}_{q_r}[\pi]=\alpha^{-1}\left(\frac{\alpha}{1+\alpha}\right)^r$, (\ref{eqn.pfC1}) equals $\gamma \left(\frac{\alpha}{1+\alpha}\right)^R$.

\newpage
\paragraph{Acknowledgements}
John Paisley and Michael I. Jordan are supported by ONR grant number N00014-11-1-0688 under the MURI program. David M. Blei is supported by ONR 175-6343, NSF CAREER 0745520, AFOSR 09NL202, the Alfred P. Sloan foundation, and a grant from Google.

\bibliographystyle{sjs}
\bibliography{nips2011bib}

\end{document}